\documentclass[11pt]{article}

\usepackage{amsxtra,amssymb,amsthm,amsmath,latexsym}
\def\R{{\mathbb R}}
\def\q{\quad}

\def\oH{\buildrel\circ\over H}
\def\oH1{\buildrel\circ\over H\kern-.02in{}^1}

\def\f{\frac}
\def\s{\sigma}
\def\b{\beta}
\def\ep{\epsilon}
\def\p{\varphi}

\def\d{{\delta}}

\begin{document}


\title{Dynamical systems method for solving operator equations
\footnote{Math subject classification: 34R30,  35R25, 35R30, 
37C35, 37L05, 37N30, 47A52, 47J06, 65M30, 65N21;
PACS 02.30.-f, 02.30.Tb, 02.30.Zz,02.60Lj, 02.60.Nm, 02.70.Pt, 05.45.-a}
}

\author{ A.G. Ramm\\
Mathematics Department, Kansas State University, \\
Manhattan, KS 66506-2602, USA\\
E:mail: ramm@math.ksu.edu     www.math.ksu.edu/\,$\widetilde{\ }$\,ramm}


\date{}
\maketitle

\begin{abstract}
 Consider an operator equation $F(u)=0$ in a real Hilbert space.
 The problem of solving this equation is ill-posed if the operator
 $F'(u)$ is not boundedly invertible, and well-posed otherwise.
 A general method, dynamical systems method (DSM) for solving linear and 
nonlinear
 ill-posed problems in a Hilbert space is presented.
 This method consists of the construction of a nonlinear
 dynamical system, that is, a Cauchy problem, which
 has the following properties:
 1) it has a global solution,
 2) this solution tends to a limit as time tends to infinity,
 3) the limit solves the original linear or non-linear
  problem.
 New convergence and discretization theorems are obtained.
 Examples of the applications of this approach are given.
 The method works for a wide range of well-posed problems as well.
\end{abstract}

\section{Introduction}

 This paper contains a recent development of the theory of DSM (dynamical 
systems method) earlier developed in papers
[2]-[12]. DSM is a general method for solving operator equations, 
especially nonlinear, ill-posed, but also well-posed operator equations. 
The author hopes that DSM will demonstrate its practical efficiency
and will allow one to solve ill-posed problems which cannot be solved by 
other methods.
This paper is intended for a broad audience: the presentation is 
simplified considerably, and is non-technical in its present form. 
Most of the results are presented  in a new way. Some of the results 
and/or proofs are new  
(Theorems 2.1, 3.1, 3.2, 4.2,  6.2, 7.1, 8.1, Remarks 4.4, 4.5, and the 
discussion of the stopping rules).
We try to emphasize the basic ideas and methods of the proofs. 

{\it What is the dynamical systems method (DSM) for solving operator 
equations?}

Consider an equation
$$
F(u):=B(u)-f=0, \quad f\in H,
\eqno{(1.1)}
$$
where $B$ is a linear or nonlinear operator in a real Hilbert space $H$.
Some of our results can be generalized to more general spaces, but these
generalizations are not discussed here. Throughout the paper we assume
that:
$$
\sup_{u\in B(u_0,R)} ||F^{(j)}(u)||\leq M_j, \quad j=1,2,
\eqno{(1.2)}
$$
where $B(u_0,R):=\{u: ||u-u_0||\leq R \}$, $F^{(j)}(u)$ is the Fr\'echet 
derivative, and 
$$
F(y)=0, \quad y\in B(u_0,R),
\eqno{(1.3)}
$$
that is, we assume existence of a solution to (1.1), not necessarily uique
globally.

Assumptions (1.2) and (1.3) are our standard assumptions below, unless 
otherwise stated. Only for well-posed problems in Section 2 we do not 
assume existence of a solution, but prove it, and sometimes we can assume
in these problems $j=1$ in (1.2), rather than $j=2$. 
In all the ill-posed problems we assume existence of the solution to 
(1.1).

Let $\dot u$ denote derivative with respect to time. Consider the 
dynamical system ( the Cauchy problem ):
$$
\dot u=\Phi(t,u), \,\,\, u(0)=u_0, 
\eqno{(1.4)}
$$
where $\Phi(t,u)$ is locally Lipschitz with respect to $u\in H$ and 
continuous 
with respect to $t\geq 0$:
$$
\sup_{u,v\in B(u_0,R), t\in [0,T]}||\Phi(t,u)-\Phi(t,v)||\leq c||u-v||, 
\quad 
c=c(R,u_0, T)>0.                     
\eqno{(1.5)}
$$
One can relax "locally Lipschitz" assumption about $\Phi$
(for example, use one-sided inequalities), but we do not 
discuss this point.
Problem (1.4) has a unique local solution if (1.5) holds. The DSM 
for solving (1.1) consists of solving (1.4), where $\Phi$ is so chosen 
that the following three conditions hold:
$$
\exists u(t) \forall t>0; \quad \exists u(\infty):=\lim_{t\to \infty}u(t); 
\quad F(u(\infty))=0.
\eqno{(1.6)}
$$
 Some of the basic results of this paper are the Theorems which provide 
the choices of $\Phi$ for which (1.6) holds, and 
the technical tools (Theorems 4.1 and 6.1) basic for our proofs.

Problem (1.1) with noisy data
$f_\delta$, $||f_\delta -f||\leq \delta$, given in place of $f$, generates
the problem:
$$
\dot u_\d=\Phi_\d(t,u_\d), \,\,\, u_\d(0)=u_0,
\eqno{(1.7)}
$$
The solution $u_\d$ to (1.7), calculated at $t=t_\d,$ will have the 
property
$$
\lim_{\d \to 0}||u_\d(t_\d)-y||=0.
\eqno{(1.8)}
$$
 The choice of $t_\d$ with this property is called the stopping rule.
One has usually $\lim_{\d \to 0}t_\d=\infty$.

In Section 2 we discuss well-posed problems (1.1), 
that is, the problems for which
$$
\sup_{u\in B(u_0,R)}||[F'(u)]^{-1}||\leq m_1,
\eqno{(1.9)}
$$
and in the other 
sections  ill-posed problems (1.1), for which (1.9) fails, are 
discussed.

The motivations for this work are: 

1) to develop a general method for 
solving operator equations, especially nonlinear and ill-posed, 

and

2) to develop a general approach to constructing convergent iterative 
schemes for solving these equations.

If (1.6) holds, and if one constructs a convergent discretization scheme 
for solving Cauchy problem (1.4), then one gets a convergent iterative 
scheme for solving the original equation (1.1).

 \section{Well-posed problems}
Consider (1.1), let (1.2) hold, and assume
$$
(F'(u) \Phi (t,u), F(u)) \leq -g_1(t)||F(u)||^a \quad \forall u\in 
B(u_0,R), \quad \int_0^\infty g_1dt=\infty,
\eqno{(2.1)}
$$
 where $g_1>0$ is an integrable function, $a>0$ is a constant.
Assume
$$
||\Phi (t,u)||\leq g_2(t)|| F(u)||, \quad \forall u\in B(u_0,R),
\eqno{(2.2)}
$$
where $g_2>0$ is such that
$$
G(t):=g_2(t)\exp(-\int_0^tg_1ds) \in L^1(\R_+).
\eqno{(2.3)}
$$
{\bf Remark: } Sometimes the assumption (2.2) can be used
in the following modified form:
$$
||\Phi (t,u)|| \leq g_2(t)||F(u)||^b \quad \forall u\in B,
\eqno{(2.2')}
$$
where $b>0$ is a constant. The statement and proof of Theorem 2.1 can be
easily
adjusted to this assumption.

Our first basic result is the folowing:

{\bf Theorem 2.1.} {\it i) If (2.1)-(2.3) hold, and 
$$ ||F(u_0)|| \int_0^\infty G(t)dt\leq R, \quad a=2,
\eqno{(2.4)}$$
then (1.4) has a global solution, (1.6) holds, (1.1) has a solution
$y=u(\infty)\in B(u_0,R)$, and 
$$ 
||u(t)-y||\leq ||F(u_0)||\int_t^\infty G(x)dx,\quad ||F(u(t))||\leq 
||F(u_0)||\exp(-\int_0^t g_1(x)dx).
\eqno{(2.5)}$$

ii) If (2.1)-(2.3) hold, $0<a<2$, and
$$ 
||F(u_0)||\int_0^Tg_2ds\leq R,
\eqno{(2.6)}$$
where $T>0$ is defined by the equation
$$
\int_0^T g_1(s)ds=||F(u_0)||^{2-a}/(2-a),
\eqno{(2.7)}$$
then (1.4) has a global solution, (1.6) holds, (1.1) has a solution
$y=u(\infty)\in B(u_0,R)$, and $u(t)=y$ for $t\geq T$.

iii) If (2.1)-(2.3) hold, $a>2$, and
$$
\int_0^\infty g_2(s)h(s)ds\leq R,
\eqno{(2.8)}$$
where 
$$ [||F(u_0)||^{2-a} +(a-2)\int_0^tg_1(s)ds]^{\frac 1
{2-a}}:=h(t), \quad \lim_{ t\to \infty} h(t)= 0,
\eqno{(2.9)}$$
then (1.4) has a global solution, (1.6) holds, (1.1) has a solution
$y=u(\infty)\in B(u_0,R)$, and 
$$ 
||u(t)-u(\infty)||\leq \int_t^\infty g_2(s)h(s)ds\to 0 
\eqno{(2.10)}$$
as $t\to \infty$.
}

Let us sketch the proof.

{\bf Proof of Theorem 2.1.} The assumptions about  $\Phi$ imply local 
existence and uniqueness of the
solution $u(t)$ to (1.4). To prove global existence of
$u$, it is sufficient to prove a uniform with respect to $t$ bound on   
$||u(t)||$. Indeed, if the maximal interval of the existence of $u(t)$
is finite, say $[0,T)$, and  $\Phi(t,u)$ is locally Lipschitz with respect 
to $u$, then $||u(t)|| \to \infty$ as $t\to T$.
  
Assume $a=2$. Let $g(t):=||F(u(t))||$. Since $H$ is real, one uses (1.4)  
and (2.1)
to get $g\dot g=(F'(u)\dot u, F)\leq -g_1(t)g^2$, so $\dot g\leq
-g_1(t)g$, and integrating one gets the second inequality (2.5), because 
$g(0)=||F(u_0)||$.
Using (2.2), (1.4) and  the second inequality (2.5), one gets:
$$ ||u(t)-u(s)||\leq g(0)\int_s^t G(x)dx, \quad G(x):=g_2(x)\exp(-\int_0^x
g_1(z)dz).
\eqno{(2.5')}
$$
Because $G\in L^1(R_+)$, it follows from ($2.5'$) that the limit
$y:=\lim_{t\to \infty} u(t)=u(\infty)$ exists, and $y\in B$ by (2.4).
From  the second inequality (2.5) and the continuity of $F$ one gets 
$F(y)=0$, so $y$ solves
(1.1).
Taking $t\to \infty$ and setting $s=t$ in ($2.5'$) yields the first 
inequality (2.5). The inclusion $u(t)\in B$ for all $t\geq 0$ follows from
(2.4) and ($2.5'$).
The first part of Theorem 2.1 is proved. The proof of the other parts is 
similar. $\Box$

There are many applications of this theorem. We mention just a few,
and assume that $g_1=c_1=const>0$ and $g_2=c_2=const>0$.

 {\bf Example 1. Continuous Newton-type method (Gavurin, (1958)):}

$\Phi= -[F'(u)]^{-1} F(u)$. Assume that (1.9) holds, then $c_1=1, 
c_2=m_1$, (2.4) takes the form (*) $ m_1(R)||F(u_0)||\leq R$, and  
(*) implies that (1.4) has a global solution, (1.6) and (2.5) hold, 
and (1.1) 
has a solution in $B(u_0,R)$. 

{\bf Example 2. Continuous simple iterations method:}

Let $\Phi=-F$, and assume $F'(u)\geq c_1(R)>0$ for all $u\in B(u_0,R)$.
 Then $c_2=1$, $c_1=c_1(R)$,  (2.4) is: 
$  [c_1(R)]^{-1}||F(u_0)||\leq R$, and the conclusions of Example 1 hold.

{\bf Example 3. Continuous gradient method:}

Let $\Phi=-[F']^{*}F$, (1.2) and (1.9) hold, 
 $c_1=m_1^{-2}$, $c_2=M_1(R)$, (2.4) is (**) \\
$M_1m_1^2||F(u_0)||\leq R$,
and (**) implies the conclusions of Example 1.

{\bf Example 4. Continuous Gauss-Newton method:}
Let $\Phi=-([F']^{*}F')^{-1}[F']^{*}F$, (1.2) and (1.9) hold,
$c_1=1$, $c_2=m_1^2M_1$, (2.4) is (***)
$  M_1m_1^2 ||F(u_0)||\leq R$, and (***) implies the conclusions
of Example 1.

{\bf Example 5. Continuous modified Newton method:}
Let $\Phi= -[F'(u_0)]^{-1} F(u)$. Assume $||[F'(u_0)]^{-1}||\leq m_0$,
and let (1.2) hold. Then $c_2=m_0$. Choose $R=(2M_2m_0)^{-1}$, and 
$c_1=0.5$. Then (2.4) is $2m_0||F(u_0)||\leq (2M_2m_0)^{-1}$, that is,
$ 4 m_0^2 M_2||F(u_0)||\leq 1$. Thus, if $4 m_0^2 M_2||F(u_0)||\leq 1$,
then the conclusions of Example 1 hold.

{\bf Example 6. Descent methods.}

Let $\Phi=-\frac{f}{(f',h)}h$, where $f=f(u(t))$ is a
differentiable functional  $f: H\to [0,\infty)$, and $h$ is an
element of $H$.
From (1.4) one gets $\dot f=(f',\dot u)=-f$. Thus $f=f_0e^{-t}$, where 
$f_0:=f(u_0)$.
Assume $||\Phi||\leq c_2|f|^b,\, b>0$. Then $||\dot u||\leq
c_2|f_0|^be^{-bt}$. Therefore $u(\infty)$ does exist, $f(u(\infty))=0$,
and $||u(\infty)-u(t)||\leq ce^{-bt}$, $c=const>0$.

If $h=f'$, and $f=||F(u)||^2$, then $f'(u)=2[F']^*(u)F(u)$,   
$\Phi=-\frac{f}{||f'||^2}f'$, and (1.4) is a descent
method. For this $\Phi$ one has $c_1=\frac 1 2$, and $c_2=\frac {m_1}
2$, where $m_1$ is defined in (1.9). Condition (2.4) is:
$  m_1 ||F(u_0)||\leq R.$ If this inequality holds, then
the conclusions of Example 1 hold.

In Example 6 we have obtained some results from [1].
Our approach is more general
than the one in [1], since the choices of $f$ and $h$ do not allow one,  
for example, to obtain $\Phi$ used in Example 5.

\section{Linear ill-posed problems}

We assume that (1.9) fails. Consider
$$
Au=f.
\eqno{(3.1)}$$
Let us denote by A) the folowing assumption:
 
{\it A):}  {\it $A$ is a linear, bounded operator in 
$H$, defined on all of $H$,
the range $R(A)$ is not closed, so (3.1) is an ill-posed 
problem, there is a $y$ such that $Ay=f$, $y\perp N$, where $N$ is the 
null-space of $A$. }

Let $B:=A^*A$, 
$q:=By=A^*f$, $A^*$ is the adjoint of $A$.
Every solution to (3.1) solves
$$
Bu=q,
\eqno{(3.2)}$$
and, if $f=Ay$, then every  solution to (3.2) solves (3.1). Choose a 
continuous, 
monotonically
decaying to zero function  $\ep(t)>0$, on $\R_+$.

Sometimes it is convenient to assume that
$$
\lim_{t\to \infty} (\dot\ep \ep^{-2})=0.
\eqno{(3.3)}$$
For example, the functions $\ep=c_1(c_0+t)^{-b}, \, 0<b<1,$ where $c_0$ 
and $c_1$ are 
positive constants, satisfy (3.3).
There are many such functions. One
can prove ([4], [8]) the following:

{\bf Claim:} {\it If $\ep(t)>0$ is a continuous monotonically decaying 
function on 
$\R_+$, $\lim_{t\to \infty}\ep(t)=0$, and (3.3) holds, then
$$
\int_0^\infty \ep ds=\infty.
\eqno{(3.3')}
$$ 
}
In this Section we
do not use assumption (3.3): in the proof of Theorem 3.1
one uses only the monotonicity of a continuous function $\ep>0$ and 
$(3.3')$.
One can drop assumption $(3.3')$, but then convergence is proved in 
Theorem 3.1
to some element of $N$, not necessarily to the normal solution $y$, that 
is, to the solution orthogonal to $N$, or, which is the same, to the 
minimal norm solution to (3.1). However, (3.3) is used (in a slightly 
weaker form) in Section 4. 

Consider problems (1.4) and (1.7) with 
$$\Phi:=-[Bu+\ep(t)u-q], \quad \Phi_\d= -[Bu_\d+\ep(t)u_\d-q_\d],
\eqno{(3.4)}
$$
 where
$||q-q_\d||\leq ||A^*||\d:=C\d$. Without loss of generality one may assume 
this $C=1$, which we do in what follows. Our main result in Sec. 3, is
Theorem 3.1, stated below. It yields the following: 

{\bf Conclusion: }{\it Given noisy data $f_\d$, every linear ill-posed 
problem 
(3.1)
under the assumptions A)  can be stably solved by the DSM.}

The result presented in Theorem 3.1 is essentially obtained in [8], but
our proof is different and much shorter.

{\bf Theorem 3.1.} {\it Problem (1.4) with $\Phi$ from (3.4) has
a unique global solution $u(t)$, (1.6) holds, and $u(\infty)=y$.
Problem (1.7) with $\Phi_\d$ from (3.4), has a unique global 
solution $u_\d(t)$, and there exists $t_\d$, such that
$$
\lim_{\d \to 0}||u_\d(t_\d)-y||=0.
\eqno{(3.5)}
$$
 This $t_\d$ can be chosen, for example, as a root of the equation
$$
\ep(t)=\d^{b}, \quad b\in (0,1),
\eqno{(3.6)}
$$
or of the eqation (3.6'), see below.
}

{\bf Proof of Theorem 3.1.}  Linear equations (1.4) with bounded operators
have unique global solutions. If $\Phi=-[Bu+\ep(t)u-q]$, then
the solution $u$ to (1.4) is 
$$
u(t)=h^{-1}(t) U(t)u_0+h^{-1}(t)\int_0^{||B||}\exp(-t\lambda) 
\int_0^t e^{s\lambda}h(s)ds \lambda dE_{\lambda}y,
\eqno{(3.7)}
$$
where $h(t):=\exp (\int_0^t \ep(s)ds)\to \infty$ as $t \to \infty$,
$E_{\lambda}$ is the resolution of the identity corresponding to the 
selfadjoint operator $B$, and $U(t):=e^{-tB}$ is a nonexpansive operator,
because $B\geq 0$. Actually, (3.7) can be used also when $B$ is unbounded,
$||B||=\infty$.

Using L'Hospital's rule one checks that
$$\lim_{t \to \infty} \frac {\lambda \int_0^t
e^{s\lambda}h(s)ds}{e^{t\lambda}h(t)}=\lim_{t \to \infty} \frac {\lambda
e^{t\lambda}h(t)}
{\lambda e^{t\lambda}h(t)+e^{t\lambda}h(t)\epsilon (t)}=
1 \quad \forall \lambda >0,
\eqno{(3.8)}$$
provided only that $\ep(t)>0$ and $\lim_{t\to \infty}\ep(t)=0$.
From (3.7), (3.8), and the Lebesgue dominated convergence theorem, one 
gets
$u(\infty)= y-Py$, where $P$ is the orthogonal projection operator on 
the null-space of $B$. Under our assumptions
A), $P=0$, so $u(\infty)=y$.
If $v(t):=||u(t)-y||$, then $\lim_{t\to \infty}v(t)=0$. 
In general, the rate of convergence of $v$ to zero can be 
arbitrarily slow for a suitably chosen $f$. Under an additional a priori
assumption on $f$ (for example, the source type assumptions), this rate 
can be 
estimated.

Let us describe a method for deriving a stopping rule.
One has:
$$
||u_\d(t)-y||\leq ||u_\d(t)-u(t)||+v(t).
\eqno{(3.9)}
$$
Since $\lim_{t\to\infty}v(t)=0$, any choice of $t_\d$ such that
$$
\lim_{t_\d\to\infty}||u_\d(t_\d)-u(t_\d)||=0,
\eqno{(3.10)}
$$
gives a stopping rule: for such $t_\d$ one has
$\lim_{\d \to 0}||u_\d(t)-y||=0$.

To prove that (3.6) gives such a rule, it is sufficient to check that
$$
||u_\d(t)-u(t)||\leq \f \d {\ep(t)}.
\eqno{(3.11)}
$$
Let us prove (3.11). Denote $w:=u_\d -u$. Then 
$$
\dot w=-[Bw+\ep w -p], \quad w(0)=0,\quad ||p||\leq \d.
\eqno{(3.12)}
$$
Integrating (3.12), and using the property $B\geq 0$, one gets (3.11).

Alternatively, multiply (3.12) by $w$, let $||w||:=g$, use $B\geq 0$, and 
get $\dot g\leq -\ep(t) g +\d,\, g(0)=0.$ Thus, $g(t)\leq \d 
\exp(-\int_0^t\ep ds)\int_0^t \exp(\int_0^s\ep d\tau)ds\leq \f \d 
{\ep(t)}.$ A more precise estimate, also used at the end of the proof of 
Theorem 3.2 below, yields:
$$
||u_\d(t)-u(t)||\leq \f \d {2\sqrt{\ep(t)}},
\eqno{(3.11')}
$$
and the corresponding stopping time $t_\d$ can be taken as the root
of the equation:
$$
2 \sqrt{\ep(t)}=\d^{b}, \quad b\in (0,1).
\eqno{(3.6')}
$$
 Theorem 3.1 is proved. $\Box$.

If the rate of decay of $v$ is known, then a more efficient stopping rule 
can be derived: $t_\d$ is the minimizer of the problem:
$$
v(t)+\d [\ep(t)]^{-1}=\min.
\eqno{(3.13)}
$$

For example, 
if $v(t)\leq c\ep^a(t)$, then $t_\d$ is the 
root of the equation $\ep(t)=(\frac {\d}{ca})^{\frac 1 {1+a}}$,
which one gets from (3.13) with $v=c\ep^a$. 
 
One can also use a stopping rule based on an a posteriori choice of
the stopping time, for example, the choice by a discrepancy principle.

{\it A method, much more efficient numerically than Theorem 3.1, is given 
below 
in Theorem 4.2.} 

For linear equation (3.2) with exact data this method 
uses
(1.4) with
$$
\Phi= -(B+\ep(t))^{-1}[Bu+\ep(t)u -q]=-u+(B+\ep(t))^{-1}q,
\eqno{(3.14)}
$$
and for noisy data it uses (1.4) with $\Phi_\d=-u_\d+(B+\ep(t))^{-1}q_\d$.
The linear operator $B\geq 0$ is monotone, so Theorem 4.2 is applicable.
For exact data (1.4) with $\Phi$, defined in (3.14), yields:
$$
\dot u=-u+(B+\ep(t))^{-1}q,\quad u(0)=u_0,
\eqno{(3.15)}
$$
and (1.6) holds if $\ep(t)>0$ is monotone, continuous, decreasing to $0$
as $t\to \infty$.  

Let us formulate the result:

{\bf Theorem 3.2.} {\it Assume A), and let $B:=A^*A$, $q:=A^*f$. Assume
$\ep(t)>0$ to be a continuous, monotonically decaying to zero function
on $[0,\infty)$. Then, for any $u_0\in H$, problem (3.15) has a 
unique global solution, $\exists u(\infty)=y$, $Ay=f$, and $y$ is the 
minimal-norm solution to (3.1). If $f_\d$ is given in place of $f$,
$||f-f_\d||\leq \d$, then (3.5) holds, with $u_\d(t)$ solving
(3.15) with $q$ replaced by $q_\d:=A^*f_\d$,  and  $t_\d$ is chosen, for 
example, as the root of ($3.6'$) (or by a discrepancy principle).
}

{\bf Proof of Theorem 3.2.}  One has $q=Bz$, where $Az=f$, and the 
solution to (3.15) is
$$
u(t)=u_0e^{-t}+e^{-t}\int_0^te^{s}(B+\ep(s))^{-1}Bzds:=u_0e^{-t}+
\int_0^{||B||}j(\lambda,t)dE_\lambda z
\eqno{(3.16)}
$$
where
$$
j(\lambda,t):=\int_0^t\f {\lambda e^{s}}{[\lambda +\ep(s)]e^t}ds,
\eqno{(3.17)}
$$
and $E_\lambda$ is the resolution of the identity of the selfadjoint 
operator $B$. One has
$$
0\leq j(\lambda, t)\leq 1, \q \lim_{t\to \infty}j(\lambda, t)=1 \,\,\, 
\lambda>0, \q j(0, t)=0.
\eqno{(3.18)}
$$
 From (3.16)-(3.18) it follows that $\exists u(\infty)$, 
$u(\infty)=z-P_Nz=y$, where $y$ is the minimal-norm solution to (3.1),
$N:=N(B)=N(A)$ is the null-space of $B$ and of $A$, and $P_N$ is the 
orthoprojector onto $N$ in $H$. This proves the first part of Theorem 3.2.

To prove the second part, denote $w:= u_\d-u$, $g:=f_\d-f$, where we 
dropped 
the dependence on $\d$ in $w$ and $g$ for brevity. Then 
$\dot w =-w +(B+\ep(t))^{-1}A^*g, \, w(0)=0.$ Thus
$w=e^{-t}\int_0^te^s (B+\ep(s))^{-1}A^*gds$, so $||w||\leq 
\d \,e^{-t}\int_0^t\f{e^s}{2\sqrt{\ep(s)}}ds\leq \f \d {2\sqrt{\ep(t)}}$,
where the known estimate (see e.g. [5]) was used: $||(B+\ep)^{-1}A^*||\leq 
\f 1 {2\sqrt{\ep}}$. Theorem 3.2 is proved. $\Box$

\section{Nonlinear ill-posed problems with monotone operators}

There is a large literature on the equations (1.1) and (1.4) with
monotone operators. In the result we present
the problem is nonlinear and ill-posed, the new technical tool, Theorem 
4.1, is used, and the stopping rules are discussed. 

Consider
(1.4) with monotone $F$ under standard assumptions (1.2) and (1.3),
and 
$$\Phi=-A_{\ep(t)}^{-1}(u)\bigl[F(u(t))+\ep(t)(u(t)-\tilde u_0)\bigr],
\eqno{(4.1)}
$$
where $A=A(u):=F'(u)$, $A^*$ is its adjoint, $\ep(t)$ is the same as in 
Theorem 3.2, and in Theorem 4.2 $\ep(t)$ is further specified, $\tilde 
u_0\in B(u_0,R)$ is an element we can choose to improve 
the numerical performance of the method.
If noisy data are given, then, as in Sec.3, we take 
$$F(u):=B(u)-f, \quad 
\Phi_\d=-A_{\ep(t)}^{-1}(u_\d)\bigl[B(u_\d(t))-f_\d+\ep(t)(u_\d(t)-\tilde 
u_0)\bigr],
$$
 where $||f_\d -f||\leq \d,$ 
$B$ is a monotone nonlinear operator, $B(y)=f$, and $u_\d$ solves (1.7).

To prove that (1.4) with the above $\Phi$ has a global solution and (1.6)
holds, we use the following:

{\bf Theorem 4.1.} {\it Let $\gamma(t), \s(t), \b(t)\in C[t_0,\infty)$ for 
some real number $t_0$.
If there exists a positive  function $\mu(t)\in C^1[t_0,\infty)$ such that
$$
0\le\s(t)\le\f{\mu(t)}{2}[\gamma(t)-\frac{\dot{\mu}(t)}{\mu(t)}],\quad
\b(t)\le\frac{1}{2\mu(t)}[\gamma(t)-\frac{\dot{\mu}(t)}{\mu(t)}],\quad
g_0\mu(t_0)<1,
\eqno{(4.2)}
$$
where $g_0$ is the initial condition in (4.3), then a nonnegative 
solution $g$ to the following differential inequality:
$$
\dot{g}(t)\le -\gamma(t) g(t)+\s(t)g^2(t)+\b(t),\quad 
g(t_0)=g_0,
\eqno{(4.3)}
$$
satisfies the estimate:
$$
0\leq g(t)\,\le\, \f{1-\nu(t)}{\mu(t)}\,<\, \f{1}{\mu(t)},
\eqno{(4.4)}
$$
for all $t\in [t_0,\infty)$, where
$$
0<\nu(t)=\left(\f{1}{1-\mu(t_0)g(t_0)}+
\f{1}{2}\int_{t_0}^t\left(\gamma(s)-\f{\dot{\mu}(s)}{\mu(s)}\right)ds
\right)^{-1}.
\eqno{(4.5)}
$$
}
There are several novel features in this result. First, differential 
equation, which one gets from (4.3) by replacing the inequality sign
by the equality sign, is a Riccati equation, whose solution may blow up in 
a finite time, in general. Conditions (4.2) guarantee the 
global existence of the solution to this Riccati equation with the initial 
condition (4.3). Secondly, this Riccati differential equation 
cannot be integrated analytically by separation of variables.
Thirdly, the coefficient $\sigma(t)$ may grow to infinity as
$t\to \infty$, so that the quadratic term does not necessarily has a
small coefficient, or the coefficient smaller than $\gamma(t)$.
Without loss of generality one may assume $\beta(t)\geq 0$ in Theorem 4.1.
The proof of Theorem 4.1 is given in [3].

The main  result of this Section is new. It claims a global convergence
in the sense that no assumptions on the choice of the initial
approximation $u_0$ are made. Usually one assumes that $u_0$ is 
sufficiently close to the solution of (1.1) in order to prove convergence.
We take $\tilde u_0=0$ in Theorem 4.2, because in this theorem 
$\tilde u_0$ does not play any role. The proof is valid for any choice of
$\tilde u_0$, but then the definition of $r$ in Theorem 4.2 is changed. 

{\bf Theorem 4.2.} {\it If (1.2) and (1.3) hold, $\tilde u_0=0$,
$R=3r$, where $r:=||y||+||u_0||$,
and $y\in N:=\{z: F(z)=0\}$ is the (unique) minimal norm solution
to (1.1), then, for any choice of $u_0$,  
problem (1.4) with $\Phi$ defined in (4.1), $\tilde u_0=0$,
and $\ep(t)=c_1(c_0+t)^{-b}$ with some positive constants $c_1, c_0,$ and  
$b$,  specified in the proof of Theorem 4.2,
has a global solution, this solution stays in the ball
$B(u_0, R)$ and (1.6) holds.  
If $u_\d(t)$ solves (1.4) with $\Phi_\d$ in place of $\Phi$, 
then there is a $t_\d$ such that $\lim_{\d \to 0}||u_\d(t_\d)-y||=0$.
}  

{\bf Proof of Theorem 4.2.} Let us sketch the steps of the proof.
Let $V$ solve the equation
$$
F(V)+\ep(t) V=0.
\eqno{(4.6)}
$$
Under our assumptions on $F$, it is 
well known that: i) (4.6) has a unique solution
for every $t>0$, and ii) $\sup_{t\geq 0}||V||\leq ||y||$, (cf [3]). 
If $F$ is Fr\'echet 
differentiable, then $V$ is differentiable, and $||\dot V(t)||\leq 
||y|| |\dot \ep(t)|/\ep(t)$. It is also known that if (1.3) holds, then 
$\lim_{t \to \infty}||V(t)-y||=0.$
We will show that the global solution $u$ to (1.4),
with the $\Phi$ from (4.1), does exist, and $\lim_{t \to 
\infty}||u(t)-V(t)||=0$. This is done by deriving a differential 
inequality for $w:=u-V$, and by applying Theorem 4.1 to $g=||w||$.
Since $||u(t)-y||\leq ||u(t)-V(t)||+||V(t)-y||$, it then follows that
(1.6) holds. We also check that $u(t)\in B(u_0, R)$, where 
$R:=3(||y||+||u_0||)$, for any choice of $u_0$ and a suitable choice of 
$\ep$.

Let us derive the differential inequality for $w$. One has
$$\dot w=-\dot V -A_{\ep(t)}^{-1}(u)\bigl[F(u(t))-F(V(t))+\ep(t)w],
\eqno{(4.7)}
$$
and $F(u)-F(V)=Aw+K$, where $||K||\leq M_2g^2/2$, $g:=||w||$
and $M_2$ is the constant from (1.2).
Multiply (4.7) by $w$, use the monotonicity of $F$, that is,
the property $A\geq 0$,
and the estimate $||\dot V ||\leq ||y|| |\dot \ep|/\ep$, and get:
$$
\dot g\leq -g +\frac {0.5 M g^2}{\ep} +||y||\frac{|\dot \ep|}{\ep}, 
\eqno{(4.8)}
$$
where $M:=M_2$. 
Inequality (4.8) is of the type (4.3): $\gamma=1$, $\s=0.5 M/\ep$,
$\b=||y||\frac{|\dot \ep|}{\ep}$. Choose 
$$
\mu(t)=\f {2M} {\ep(t)}.
\eqno{(4.9)}
$$
 Clearly $\mu\to \infty$ as $t\to \infty$.
Let us check three conditions (4.2). One has 
$\frac{\dot{\mu}(t)}{\mu(t)}=\f{|\dot \ep|}{\ep}$. Take $\ep=
c_1(c_0+t)^{-b}$, where $c_j>0$ are constants, $0<b<1$,
and choose these constants so that $\f {|\dot \ep|}{\ep} <\f 1 2,$
for example, $\frac b{c_0}=\frac 1 4$.  
Then the first condition (4.2) is satisfied.
The second condition (4.2)
holds if  
$$
8M||y|| |\dot \ep|\ep^{-2}\leq 1.
\eqno{(4.10)}
$$ 
One has $\ep(0)=c_1c_0^{-b}$. Choose 
$$\ep(0)=4Mr.
\eqno{(4.11)}
$$
Then
$$ |\dot \ep|\ep^{-2}=bc_1^{-1}(c_0+t)^{b-1}\leq 
bc_0^{-1}c_1^{-1}c_0^b=\frac 1 {4\ep(0)}=\frac 1 {16Mr},
\eqno{(4.12)}
$$ 
so (4.10) holds.
Thus, the second condition (4.2) holds.
The last
condition (4.2) holds because  
$$
\f {2M||u_0-V_0||} {\ep(0)}\leq \frac {2Mr}{4Mr}=\frac 1 2<1.
$$
By Theorem 4.1 one concludes that 
$g=||w(t)||<\frac {\ep(t)}{2M}\to 0$ when $t\to \infty$, and
$$
||u(t)-u_0||\leq g+||V-u_0||\leq g(0)+r\leq 3r.
\eqno{(4.13)}
$$ 
This estimate implies the global 
existence of the solution to (1.4), because if $u(t)$ would have a finite
maximal interval of existence, $[0,T)$, then $u(t)$ could not stay bounded 
when $t\to T$, which contradicts the boundedness of $||u(t)||$,
and from (4.13) it follows that $||u(t)||\leq 4r$.
We have proved the first part of Theorem 4.2, namely properties (1.6).  
$\Box$

To derive a stopping rule we argue as in Sec.3. One has:
$$
||u_\d(t)-y||\leq ||u_\d(t)-V(t)||+||V(t)-y||.
\eqno{(4.14)}
$$
We have already proved that $\lim_{t\to \infty}v (t):=
\lim_{t\to \infty}||V(t)-y||= 0$.
The rate of decay of $v$ can be arbitrarily slow, in general.
Additional assumptions, for example, the source-type ones, can be used to 
estimate the rate of decay of $v(t)$. One derives differential inequality 
(4.3) for $g_\d:=||u_\d(t)-V(t)||$, and estimates 
$g_\d$ using (4.4). The analog of (4.8) for $g_\d$ contains additional
term $\f \d \ep$ on the right-hand side. If $\f \d {\ep^2}\leq \f 1 {16M},$ 
then
conditions (4.2) hold, and $g_\d<\frac {\ep(t)} {2M}$. 
Let $t_\d$ be the root of the equation $\ep^2(t)=16M \d $.
Then $\lim_{\d\to 0}t_\d=\infty$, and (1.8) holds
because $||u_\d(t_\d)-y||\leq v(t_\d)+g_\d$, 
$\lim_{t_\d\to \infty}g_\d(t_\d)=0$ and 
$\lim_{t_\d\to \infty}v(t_\d) =0$, but the convergence in (1.8) can be 
slow. See [4] and [5] for the rate of convergence under source 
assumptions.
If the rate of decay of $v(t)$ is known, then one
chooses $t_\d$ 
as the minimizer of the problem, similar to (3.13), 
$$v(t)+g_\d(t)=min,
\eqno{(4.15)}
$$  where the minimum is 
taken over $t>0$ for a fixed small $\d>0$. This yields a quasioptimal 
stopping rule.
Theorem 4.2 is proved. $\Box$ 

In [4] a local convergence result, similar to the first part of Theorem 
4.1, was obtained, that is, $||u_0-y||$ was assumed sufficiently small, 
and no discussion of noisy data was given.

Let us give another result:

{\bf Theorem 4.3.} {\it Assume that $\Phi=-F(u)-\ep(t)u$, $F$ is monotone,
$\ep(t)$ as in Theorem 3.2, and (3.3), (1.2) and (1.3) hold. Then
(1.6) holds.
}

{\bf Proof of Theorem 4.3.} As in the proof of Theorem 4.2, it is 
sufficient to prove that $\lim_{t\to \infty}g(t)=0$, where $g, w,$ and $V$
are the same as in Theorem 4.2, and $u$ solves (1.4) with the $\Phi$
defined in Theorem 4.3. Similarly to the derivation of (4.7), one gets:
$$
\dot w=-\dot V- [F(u)-F(V)+\ep(t) w]. 
\eqno{(4.16)}
$$
Multiply (4.16) by $w$, use the monotonicity of $F$, the estimate
$||\dot V||\leq \f {|\dot \ep(t)|}{\ep (t)}||y||$, which was used also in 
the proof of Theorem 4.2, and get:
$$
\dot g\leq -\ep(t) g+\f {|\dot \ep(t)|}{\ep (t)}||y||.
\eqno{(4.17)}
$$
This implies 
$$
 g(t)\leq e^{ -\int_0^t\ep(s)ds}[ g(0)+\int_0^te^{\int_0^s\ep(x)dx}\f 
{|\dot \ep(s)|}{\ep (s)}||y||ds].
\eqno{(4.18)}
$$
From our assumptions relation ($3.3'$) follows, and (4.18) together with 
(3.3)
and ($3.3'$) imply $\lim_{t\to \infty}g(t)=0$. Theorem 4.3 is proved. 
$\Box$

{\bf Remark 4.4.} One can drop assumption (1.2) in Theorem 4.3 and assume
only that $F$ is a monotone hemicontinuous operator defined on all of $H$.

{\bf Claim 4.5:} {\it If $\ep(t)=\ep=const>0$, then $\lim_{\ep\to 
0}||u(t_\ep)-y||=0$, where $u(t)$ solves (1.4) with $\Phi:= -F(u)-\ep u$,
and $t_\ep$ is any number such that $\lim_{\ep \to 0}\ep t_\ep=\infty$.}

{\bf Proof of the claim.} One has $||u(t)-y||\leq 
||u(t)-V_\ep||+||V_\ep-y||$, where $V_\ep$ solves (4.6) with
$\ep(t)=\ep=const>0$. Under our  
assumptions on $F$, equation (4.6) has a unique solution, and 
$\lim_{\ep \to 0}||V_\ep-y||=0$. So, to prove the claim, it is sufficient 
to prove that $\lim_{\ep \to 0}||u(t_\ep)-V_\ep||=0$, provided that
$\lim_{\ep \to 0}\ep t_\ep=\infty$.
Let $g:=||u(t)-V_\ep||$, and $w:=u(t)-V_\ep$. Because $\dot V_\ep=0$, one
has the equation: $\dot w=-[F(u)-F(V_\ep)+\ep w]$.
Multiplying this equation by 
$w$, and using the monotonicity of $F$, one gets
$\dot g\leq -\ep g$, so $g(t)\leq g(0)e^{-\ep t}$. Therefore
$\lim_{\ep \to 0}g(t_\ep)=0$, provided that
$\lim_{\ep \to 0}\ep t_\ep=\infty$. The claim is proved. $\Box$

{\bf Remark 4.6.} One can prove claims i) and ii), formulated below
formula (4.6), using DSM version
presented in Theorem 8.1 below.

{\bf Claim 4.7:} {\it  Assume that $F$ is monotone, (1.2) holds,
and $F(y)=0$. Then claims i) and ii), formulated below
formula (4.6), hold.}

{\bf Proof. }  First, note that ii) follows from i)
easily, because the assumptions $F(y)=0$, $F$ is monotone, and $\ep>0$, 
imply, after multiplying
$F(V)-F(y)+\ep V=0$ by $V-y$, the inequality $(V, V-y)\leq 0,$ from which
claim ii) follows. Claim i) follows from Theorem 8.1, proved below. $\Box$

{\bf Claim 4.8:} {\it  Assume that the operator $F$ is monotone, 
hemicontinuous,
defined on all of $H$, $F(y)=0$, $y$ is the minimal norm element of
$N_F:=\{z: F(z)=0\}$,
$\Phi=-F(u)-\ep(t)u$, $\ep(t)>0$, is monotone, decaying to zero, and
(3.3) holds. Then (1.6) holds
for the solution to (1.4).}

{\bf Proof.} Existence of the unique global solution to
(1.4) under our assumptions is known (see e.g. [D]:
K.Deimling, Nonlinear functional analysis, Springer Verlag, Berlin,
1985, p.99). Let $w:=u-V_b$, $g:=||w||$, where $V_b$ solves
$F(V_b)+bV_b=0$, $b=const>0$. 
It is shown in the proof of Claim 4.7 that $||V_b||\leq
||y||$, and one can prove (see e.g., [2]) that  $\lim_{b\to 
0}||V_b-y||=0$. One has
$||u(t)-y||\leq ||u(t)-V_b||+||V_b-y||$. Thus, to prove
$\lim_{t\to \infty}||u(t)-y||=0$ it is sufficient to prove that
$\lim_{t\to \infty}g(t)=0$. One has $\dot w= -[F(u)+\ep(t)u
-F(V_b)-bV_b]$. Multiply this equation by $w$, use the monotonicity
of $F$ and get: $\dot g\leq -\ep(t) g +|\ep(t)-b|||y||$.
Denote $h(t):=\exp(\int_0^t\ep(s)ds).$  Then,
$$g(\xi)\leq g(0)h^{-1}(\xi)
+h^{-1}(\xi)\int_0^\xi h(s)|\ep(s)-b|ds ||y||, \quad \forall b=const\geq 
0. 
\eqno{(4.19)}$$
Clearly $\lim_{\xi\to \infty} g(0)h^{-1}(\xi)=0$, because $\lim_{t\to
\infty} h(t)=\infty$. In fact, $\ep^{-1}\leq ct +c_0$, where
$c_0:=\ep^{-1}(0)>0$, and one can choose $0<c<1$
because of (3.3), so $h\geq (ct+c_0)^{\frac 1 c}$,
and $\lim_{t\to \infty} \ep(t) h(t)= \infty$.
Choose $b=\ep (\xi)$ and
apply L'H\^ospital's rule to the  last term in (4.19).
 L'H\^ospital's rule is applicable, and
one gets:
$$\lim_{\xi\to \infty}g(\xi)=\lim_{\xi\to \infty} 
\frac {|\dot \ep(\xi)|}{\ep^2(\xi)}
\frac  {\ep(\xi)\int_0^\xi hds}{h(\xi)}=0,$$
because $\frac  {\ep(\xi)\int_0^\xi hds}{h(\xi)}$ is a bounded function
  and $\frac {|\dot \ep(\xi)|}{\ep^2(\xi)}\to 0$
as $\xi\to \infty$.
Claim 4.8  is proved. $\Box$

The result in claim 4.8 contains the result from [AR]
(Alber, Ya. and Ryasantseva, I., On regularized evolution
equations, Funct. Diff. Eqs, 7, (2000), 177-187.), where
additional assumptions are made on $\ep (t)$, global existence
of the solution to (1.4) is assumed, and the proof contains a gap,
because it is not shown that the L'H\^ospital's rule can be applied twice.

 \section{Nonlinear ill-posed problems with non-monotone operators}

Assume that $F(u):=B(u)-f$, $B$ is a non-monotone operator, 
$A:=F'(u)$, $\tilde A:=F'(y)$, $T:=A^*A$, 
$\tilde T:=\tilde A^* \tilde A$, $T_\ep:=T+\ep I$, where $I$ is 
the identity operator, $\ep$ is as in Theorem 3.2 and $\f {|\dot 
\ep(t)}{\ep(t)}<1$, 
$$
\Phi:= -T_{\ep}^{-1}(u)[A^*(B(u)-f)+\ep(u-\tilde u_0)],\q \ep=\ep(t)>0,
\eqno{(5.1)}
$$
and $\Phi_\d$ is defined similarly, with $f_\d$ replacing $f$ and $u_\d$ 
replacing $u$. 

The main result of this Section is:

{\bf Theorem 5.1.} {\it If (1.2) and (1.3) hold, $u, u_0\in B(y,R)$,
$y-\tilde u_0=\tilde Tz, \, ||z||<<1$, (that is, $||z||$ is sufficiently 
small),
and $R$ is sufficiently small, then problem (1.4) has a unique global 
solution and (1.6) holds. If $u_\d(t)$ solves (1.7),
 then there exists a $t_\d$ such that $\lim_{\d \to 0}||u_\d(t_\d)-y||=0.$
}

The derivation of the stopping rule, that is, the choice of $t_\d$,
is based on the ideas presented in Sec.4 (cf [8], [5]).

{\bf Sketch of proof of Theorem 5.1.} 
Proof of Theorem 5.1 consists of the following steps.

First we prove that $g:=||w||:=||u(t)-y||$ satisfies a differential 
inequality (4.3), and, applying (4.4), conclude that
$g(t)<\mu^{-1}(t)\to 0$ as $t\to \infty$. A new point in this derivation
(compared with the one for monotone operators) is 
the usage of the source assumption $y-u_0=\tilde Tz$.

Secondly, we derive the stopping rule using the ideas from Sec.4.
The source assumption allows one to get a rate of convergence (see [2] and 
[5]).
Details of the proof are technical and are not included.
One can see [5] for some proofs.

Let us sketch the derivation of the differential inequality for $g$.
Write $B(u)-f=B(u)-B(y)=Aw+K$, where $||K||\leq \f {M_2} 2 g^2$, and
$\ep (u-\tilde u_0)=\ep w +\ep (y-\tilde u_0)=\ep w+\ep \tilde Tz.$ Then 
(5.1) can be written as 
$$
\Phi=-w -T_\ep^{-1}A^*K -\ep T_\ep^{-1} \tilde Tz, \q \ep:=\ep(t).
\eqno{(5.2)}
$$ 
Multiplying (1.4), with $\Phi$ defined in (5.2), by $w$, one gets:
$$
g \dot g\leq -g^2 +\f {M_2} 2 ||T_{\ep(t)}^{-1}A^*||g^3 +\ep(t) 
||T_{\ep(t)}^{-1}\tilde T||||z||g.
\eqno{(5.3)}
$$
Since $g\geq 0$, one obtains:
$$
\dot g\leq -g +\f {M_2}{4\sqrt {\ep(t)}}g^2 +\ep(t) ||T_\ep^{-1}\tilde 
T||||z||,
\eqno{(5.4)}
$$
where the estimate 
$||T_\ep^{-1}A^*||\leq \f 1{2\sqrt {\ep}}$ was used.
Clearly, 
$$||T_\ep^{-1}\tilde T||\leq ||(T_\ep^{-1}-\tilde 
T_\ep^{-1})\tilde T||
+||\tilde T_\ep^{-1}\tilde T||, \quad ||\tilde T_\ep^{-1}\tilde T||\leq 1,
\quad \ep ||T_\ep^{-1}||\leq 1,
$$
and 
$$T_\ep^{-1}-\tilde T_\ep^{-1}=T_\ep^{-1}(A^*A-\tilde A^*\tilde A)\tilde 
T_\ep^{-1}.
$$
 One has:  
$$||A^*A-\tilde A^*\tilde A||\leq 2M_2 M_1 g,\quad ||z||<<1.
$$
Let $2M_1M_2||z||\leq \frac 1 2$. This is possible since $||z||<<1.$
Using the above estimates, one transforms (5.4) into the 
following inequality:
$$
\dot g\leq - \frac 1 2  g +\f {M_2}{4\sqrt {\ep(t)}}g^2 +||z|| \ep. 
\eqno{(5.5)}
$$
Now, apply Theorem 4.1 to (5.5), choosing 
$$\mu=\f {2M_2} {\sqrt {\ep}},\quad \f {|\dot \ep|}{\ep}< \f 1 2,
\quad 16M_2 ||z||\sqrt {\ep(0)}<1, \quad  and \quad \f {2M_2 
||u_0-y||}{\sqrt {\ep(0)}}<1.
$$
Then conditions (4.2) are satisfied, and Theorem 4.1 yields the estimate:
$$
g(t)<\f {\sqrt {\ep(t)}}{2M_2}.
$$ 
This is the main part of the proof of 
Theorem 5.1.  $\Box$

\section{Nonlinear ill-posed problems: avoiding inverting of  
operators in the Newton-type continuous schemes}

In the Newton-type methods for solving well-posed nonlinear problems, 
for example, in the continuous Newton method (1.4) with 
$\Phi=-[F'(u)]^{-1}F(u)$,
the difficult and expensive part of the solution is inverting the
operator $F'(u)$. In this section we give a method to avoid
inverting of this operator. This is especially important in the ill-posed 
problems, where one has to invert some regularized versions of $F'$,
and to face more difficulties than in the well-posed problems.

Consider problem (1.1) and assume (1.2), (1.3) and (1.9). Thus, we 
discuss our method in the simplest well-posed case. 

Replace (1.4) by the following Cauchy problem (dynamical system):

$$
\dot u=-QF, \quad u(0)=u_0,
\eqno{(6.1)}
$$
$$
\dot Q=-TQ+A^*, \quad Q(0)=Q_0,
\eqno{(6.2)}
$$
where $A:=F'(u)$, $T:=A^*A$, and $Q=Q(t)$ is a bounded operator in $H$.

First let us state our new technical tool: an operator version of the 
Gronwall inequality (cf [9]).

{\bf Theorem 6.1.} {\it Let
$$
\dot Q= - T(t)Q(t)+G(t),\q Q(0)=Q_0,
\eqno{(6.3)}
$$
where $T(t),$ $G(t),$ and $Q(t)$ are linear bounded operators
on a real Hilbert
space   $H$. If there exists $\ep(t)>0$ such that
$$
(T(t)h,h)\ge \ep(t)||h||^2 \q \forall h\in H,
\eqno{(6.4)}
$$
then
$$
||Q(t)||\le
e^{-\int^t_0\ep(s)ds}[||Q(0)||+\int^t_0||G(s)||e^{\int^s_0\ep(x)dx}\,ds].
\eqno{(6.5)}
$$
}

{\it Let us turn now to a proof of Theorem 6.2, formulated at the end of 
this 
Section. This theorem is the main result of Section 6.}
 
Applying (6.5) to (6.2), and using (1.2) and (1.9), which implies
$$
(T(t)h,h)\ge c ||h||^2 \q \forall h\in H, \q c=const>0,
\eqno{(6.6)}
$$
one gets:
$$
||Q(t)||\le
e^{-ct}[||Q(0)||+\int^t_0M_1e^{cs}\,ds
]\leq [||Q_0||+M_1 c^{-1}]:=c_1,
\eqno{(6.7)}
$$
as long as $u(t)\in B(u_0,R)$.

Let $u(t)-y:=w,\, ||w||:=g,\, \tilde A:=F'(y)$. Since $F(y)=0$, 
one has $F(u)=\tilde Aw + K$, where $||K||\leq 0.5 M_2 g^2:=c_0g^2$,
and $M_2$ is the constant from (1.2).
Rewrite
(6.1) as 
$$
\dot w=-Q[\tilde Aw +K].
\eqno{(6.8)}
$$
Let $\Lambda:=I-Q\tilde A$. 
Multiply (6.8) by $w$ and get 
$$g\dot g\leq -g^2 +(\Lambda w,w)+c_0g^3, \quad c_0=const>0.
\eqno{(6.9)}
$$
We prove below that
$$
\sup_{t\geq 0}||\Lambda||\leq \lambda<1.
\eqno{(6.10)}
$$
From (6.9) and (6.10) one gets the following differential inequality:
$$
\dot g\leq -\gamma g +c_0 g^2, \q 0<\gamma <1, \,\, \gamma:=1-\lambda,
\eqno{(6.11)}
$$
which implies:
$$
g(t)\leq re^{-\gamma t}, \q r:= g(0)[1-g(0)c_0]^{-1},
\eqno{(6.12)}
$$
provided that
$$
g(0)c_0<1.
\eqno{(6.13)}
$$
Inequality (6.13) holds if $u_0$ is sufficiently close to $y$.

From (6.12) and (6.11) it folows that $u(\infty)=y$. Thus, (1.6) holds.

The trajectory $u(t)\in B(u_0,R), \, \forall t>0,$ 
provided that
$$
\int_0^\infty ||\dot u||dt=\int_0^\infty ||\dot w||dt\leq r+\frac 
{c_0r^2}{2\gamma}\leq R.
\eqno{(6.14)}
$$
 This inequality holds if $u_0$ is sufficiently close to $y$, that is, $r$ 
is sufficiently small.

To complete the argument, let us prove (6.10). One has:
$$
\dot \Lambda=-\dot Q \tilde A= -T\Lambda +A^*(A-\tilde A).
\eqno{(6.15)}
$$
One has $||A-\tilde A||\leq M_2g$. Using (6.12) and Theorem 6.1, 
one gets
$$
||\Lambda||\leq e^{-ct}[ ||\Lambda_0||+rM_1M_2\int_0^t e^{(c-\gamma) 
s}ds].
\eqno{(6.16)}
$$
 Thus,
$$
||\Lambda||\leq ||\Lambda_0|| +Cr:=\lambda, \q C:=M_1M_2 \sup_{t>0}\frac 
{e^{-\gamma 
t}-e^{-ct}}{c-\gamma}.
\eqno{(6.17)}
$$
If $u_0$ is sufficiently close to $y$ and $Q_0$ is sufficiently
close to $\tilde A^{-1}$, then $\lambda>0$ can be made arbitrarily small.
 We have proved:

{\bf Theorem 6.2.} {\it If (1.2), (1.3) and (1.9) hold, $Q_0$ and $u_0$
are sufficiently close to $\tilde A^{-1}$ and $y$, respectively, then
problem (6.1)-(6.2) has a unique global solution, (1.6) holds, and $u(t)$ 
converges to $y$, which solves (1.1), exponentially fast.
}

In [9] a generalization of Theorem 6.2 is given for ill-posed problems.

\section{Iterative schemes}
 In this section we present a method for constructing convergent iterative 
schemes for a wide class of well-posed equations (1.1). Some methods for 
constructing convergent iterative
schemes for a wide class of ill-posed problems are given in [3]. There is 
an enormous literature on iterative methods.

Consider a discretization scheme for solving (1.4) with $\Phi=\Phi(u)$,
so that we assume no explicit time dependence in $\Phi$:
$$
u_{n+1}=u_n +h\Phi(u_n),\quad u_0=u_0, \quad h=const>0.
\eqno{(7.1)}
$$
One of our results from [3], concerning the 
well-posed equations (1.1) is Theorem 7.1, formulated below.
Its proof is shorter and simpler than in [3].

{\bf Theorem 7.1.} {\it Assume (1.2), (1.3), (1.9), (2.1)-(2.4) with 
$a=2$, $g_1=c_1=const>0$, $g_2=c_2=const>0$,
 $ ||\Phi'(u)||\leq L_1$, for $u\in B(y,R)$.
Then, if $h>0$ is sufficiently small, and $u_0$ is sufficiently close to 
$y$, then (7.1) produces a sequence $u_n$ for which
$$
||u_n-y||\leq Re^{-chn},\quad ||F(u_n)||\leq ||F_0||e^{-chn},
\eqno{(7.2)}
$$
where $R:=\frac {c_2||F_0||}{c_1}, \,F_0= F(u_0),\, c=const>0, \,\, 
c<c_1.$
}

{\bf Proof of Theorem 7.1.} The proof is by induction. For $n=0$ estimates 
(7.2) are 
clear. Assuming these estimates for $j\leq n$, let us prove them for 
$j=n+1$. Let $F_{n}:=F(u_n)$, and let $w_{n+1}(t)$ solve problem (1.4) on 
the interval 
$(t_n,t_{n+1}),$ $t_n:=nh$, with $w(t_n)=u_n$. By (2.5) (with 
$G=c_2e^{-c_1t}$) and (7.2) one gets:
$$
||w_{n+1}(t)-y||\leq \frac {c_2}{ c_1} ||F_n||e^{-c_1t}\leq
Re^{-cnh-c_1t}, \quad t_n\leq t \leq t_{n+1}. 
\eqno{(7.3)}
$$
One has:
$$
||u_{n+1}-y||\leq ||u_{n+1}-w_{n+1}(t_{n+1})||+  ||w_{n+1}(t_{n+1})-y||,
\eqno{(7.4)}
$$
and
$$
 ||u_{n+1}-w_{n+1}(t_{n+1})||\leq \int_{t_n}^{t_{n+1}} 
||\Phi(u_n)-\Phi(w_{n+1}(s))||ds
$$
$$\leq 
L_1c_2h\int_{t_n}^{t_{n+1}}||F(w_{n+1}(t))||dt\leq L_1c_1h^2 Re^{-cnh}, 
\eqno{(7.5)}
$$
where we have used the formula  $R:=\frac {c_2||F_0||}{c_1}$, and the 
estimate:
$$
||F(w_{n+1}(t))||\leq ||F_n||e^{-c_1(t-t_n)}\leq 
||F_0||e^{-cnh-c_1(t-t_n)}.
\eqno{(7.6)}
$$
From (7.3)-(7.6) it follows that:
$$ 
||u_{n+1}-y||\leq Re^{-cnh} (e^{-c_1h}+c_1L_1h^2)\leq Re^{-c(n+1)h},
\eqno{(7.7)}
$$
provided that
$$ 
e^{-c_1h}+c_1L_1h^2\leq e^{-ch}.
\eqno{(7.8)}
$$
Inequality (7.8) holds if $h$ is sufficiently small and $c<c_1$.
So, the first inequality (7.2), with $n+1$ in place of $n$, is proved
if  $h$ is sufficiently small and $c<c_1$.

Now 
$$
||F(u_{n+1})||\leq ||F(u_{n+1})-F(w_{n+1}(t))||+||F(w_{n+1}(t))||,\quad
t_n\leq t \leq t_{n+1}.
\eqno{(7.9)}
$$
Using (1.2) and (7.5), one gets:
$$
||F(u_{n+1})-F(w_{n+1}(t_{n+1}))||\leq M_1 
||u_{n+1}-w_{n+1}(t_{n+1})||\leq M_1c_2L_1h^2 ||F_0||e^{-cnh}.
\eqno{(7.10)}
$$
From (7.9) and (7.10) it follows that:
$$
||F(u_{n+1})||\leq ||F_0||e^{-cnh}(e^{-c_1h}+ M_1c_2L_1h^2)\leq
 ||F_0||e^{-c(n+1)h}, 
\eqno{(7.11)}
$$  provided that
$$
e^{-c_1h}+ M_1c_2L_1h^2\leq e^{-ch}.
\eqno{(7.12)}
$$
Inequality (7.12) holds if $h$ is sufficiently small and $c<c_1$.
So, the second inequality (7.2) with $n+1$ in place of $n$ is proved
if  $h$ is sufficiently small and $c<c_1$.
Theorem 7.1 is proved. $\Box$

In the well-posed case, if $F(y)=0$, the discrete Newton's method
$$
u_{n+1}=u_n-[F'(u_n)]^{-1}F(u_n),\, u_0=u(0),
$$ 
converges superexponentially if $u_0$ is sufficiently close to $y$.
Indeed,  if $v_n:=u_n-y$,
then $v_{n+1}=v_n- [F'(u_n)]^{-1}[F'(u_n)v_n+K]$ where $||K||\leq \f 
{M_2}2 ||v_n||^2.$ Thus, $g_n:=||v_n||$ satisfies the inequality:
$g_{n+1}\leq qg_n^2$, where $q:=\f {m_1M_2}2$. Therefore
$g_n\leq q^{2^n-1}g_0^{2^n},$ and if $0<qg_0<1$, then the method converges
superexponentially. 

If one uses the iterative method 
$u_{n+1}=u_n-h[F'(u_n)]^{-1}F(u_n),$ with $h\neq 1$, then,
in the well-posed case, assuming that this method converges, it converges 
exponentially, that is, slower than in the case $h=1$. 

The continuous analog of the above method  
$$\dot u= -a[F'(u)]^{-1}F(u), \quad  u(0)=u_0,
$$ 
where $a=const>0$, converges at the rate $O(e^{-at})$. Indeed, if 
$g(t):=||F(u(t))||,$
then $g \dot g=-ag^2$, so $g(t)=g_0e^{-at}$, 
$||\dot u||\leq am_1g_0e^{-at}$. Thus 
$$
||u(t)-u(\infty)||\leq m_1g_0e^{-at},\quad and \quad F(u(\infty))=0.
$$ 
In the continuous case one does not have 
superexponential convergence no matter what $a>0$ is (see [11]).

\section{ A spectral assumption }

In this section we introduce the spectral assumption which allows one
to treat some nonlinear non-monotone operators. 

{\bf Assumption S:}  {\it The set $\{r,\p: \pi -\p_0<\p<\pi +\p_0,\,\, 
\p_0>0, 
0<r<r_0\}$, where $\p_0$ and $r_0$ are arbitrarily small, fixed numbers,
consists of the regular points of the operator $A:=F'(u)$ for all $u\in 
B(u_0,R)$.}

Assumption S implies the estimate:
$$
||(F'(u)+\ep)^{-1}||\leq \f 1 {\ep \sin \p_0}, \q \ep< r_0(1-\sin\p_0),\q
\ep=const>0,
\eqno{(8.1)}
$$
because $||(A-z)^{-1}||\leq \f 1 {dist(z, s(A))}$, where $s(A)$ is the 
spectrum of a linear operator $A,$ and $dist(z, s(A))$ is the distance 
from a point $z$ of a complex plane to the spectrum. In our case, 
$z=-\ep$,  and $dist(z, s(A))=\ep \sin \p_0,$ if $\ep< r_0(1-\sin\p_0)$.

{\bf Theorem 8.1.} {\it If (1.2) and (8.1) hold, and 
$0<\ep< r_0(1-\sin\p_0)$, then problem (4.6), with $\ep(t)=\ep=const>0$,
 is solvable, problem (1.4),
with $\Phi$ defined in (4.1) and $\tilde u_0=0$, has a unique global 
solution, $\exists u(\infty)$, and $F(u(\infty))+\ep u(\infty)=0$.
Every solution to the equation $F(V)+\ep V=0$ is isolated.
}

{\bf Proof of Theorem 8.1.} Let $g=g(t):=||F(u(t))+\ep u(t)||$, where 
$u=u(t)$ solves 
locally (1.4), where $\Phi$ is defined in (4.1) and $\tilde u_0=0$.
Then:
$$
g\dot g=-((F'(u)+\ep)(F'(u)+\ep)^{-1}(F(u)+\ep u), F(u)+\ep u)=-g^2,
\eqno{(8.2)}
$$
so
$$
g=g_0e^{-t},\,\, g_0:=g(0);\q ||\dot u||\leq \f {g_0} {\ep \sin 
\p_0}e^{-t}.
\eqno{(8.3)}
$$
Thus,
$$
||u(t)-u(\infty)||\leq \f {g_0} {\ep \sin \p_0}e^{-t},\q 
||u(t)-u_0||\leq \f {g_0} {\ep \sin \p_0}, \q F(u(\infty))+\ep 
u(\infty)=0.
\eqno{(8.4)}
$$
Therefore equation
$$
F(V)+\ep V=0, \q \ep=const>0,
\eqno{(8.5)}
$$
has a solution in $B(u_0,R)$, where $R=\f {g_0} {\ep \sin \p_0}.$

Every solution to equation (8.5) is isolated. Indeed, if 
$F(W)+\ep W=0$, and $\psi:=V-W$, then $F(V)-F(W)+\ep \psi=0,$ 
so $[F'(V)+\ep]\psi +K=0$, where $||K||\leq \f {M_2}2 ||\psi||^2$.
Thus, using (8.1), one gets $||\psi||\geq \f {2\ep \sin \p_0}{M_2}$.
Consequently, if $||\psi||$ is sufficiently small, then $\psi=0$.
Theorem 8.1 is proved. $\Box$

\end{document}